\documentclass[12pt]{article}
\usepackage{amsmath,amssymb,amsthm,latexsym,amscd}

\title{Relative $e$-spectra and
relative closures \\ for families of theories\footnote{
\newline\indent \ \ \ The research is partially supported by the Grants Council (under RF President) for State Aid of
Leading Scientific Schools (grant NSh-6848.2016.1) and by
Committee of Science in Education and Science Ministry of the
Republic of Kazakhstan (Grant No. 0830/GF4).} }
\author{Sergey V.
Sudoplatov\footnote{sudoplat@math.nsc.ru}}
\date{}
\begin{document}
\maketitle

\begin{abstract}
We define the notions of relative $e$-spectra, with respect to
$E$-operators, relative closures, and relative generating sets. We
study properties connected with relative $e$-spectra and relative
generating sets.

\end{abstract}

We continue to study structural properties of combin\-a\-tions of
structures and their theories \cite{cs, cl, lut} generalizing the
notions of $e$-spectra, closures and generating sets to relative
ones. Properties of relative $e$-spectra and relative generating
sets are investigated.

\section{Preliminaries}

Throughout the paper we use the following terminology in \cite{cs,
cl}.

Let $P=(P_i)_{i\in I}$, be a family of nonempty unary predicates,
$(\mathcal{A}_i)_{i\in I}$ be a family of structures such that
$P_i$ is the universe of $\mathcal{A}_i$, $i\in I$, and the
symbols $P_i$ are disjoint with languages for the structures
$\mathcal{A}_j$, $j\in I$. The structure
$\mathcal{A}_P\rightleftharpoons\bigcup\limits_{i\in
I}\mathcal{A}_i$\index{$\mathcal{A}_P$} expanded by the predicates
$P_i$ is the {\em $P$-union}\index{$P$-union} of the structures
$\mathcal{A}_i$, and the operator mapping $(\mathcal{A}_i)_{i\in
I}$ to $\mathcal{A}_P$ is the {\em
$P$-operator}\index{$P$-operator}. The structure $\mathcal{A}_P$
is called the {\em $P$-combination}\index{$P$-combination} of the
structures $\mathcal{A}_i$ and denoted by ${\rm
Comb}_P(\mathcal{A}_i)_{i\in I}$\index{${\rm
Comb}_P(\mathcal{A}_i)_{i\in I}$} if
$\mathcal{A}_i=(\mathcal{A}_P\upharpoonright
A_i)\upharpoonright\Sigma(\mathcal{A}_i)$, $i\in I$. Structures
$\mathcal{A}'$, which are elementary equivalent to ${\rm
Comb}_P(\mathcal{A}_i)_{i\in I}$, will be also considered as
$P$-combinations.

Clearly, all structures $\mathcal{A}'\equiv {\rm
Comb}_P(\mathcal{A}_i)_{i\in I}$ are represented as unions of
their restrictions $\mathcal{A}'_i=(\mathcal{A}'\upharpoonright
P_i)\upharpoonright\Sigma(\mathcal{A}_i)$ if and only if the set
$p_\infty(x)=\{\neg P_i(x)\mid i\in I\}$ is inconsistent. If
$\mathcal{A}'\ne{\rm Comb}_P(\mathcal{A}'_i)_{i\in I}$, we write
$\mathcal{A}'={\rm Comb}_P(\mathcal{A}'_i)_{i\in
I\cup\{\infty\}}$, where
$\mathcal{A}'_\infty=\mathcal{A}'\upharpoonright
\bigcap\limits_{i\in I}\overline{P_i}$, maybe applying
Morleyzation. Moreover, we write ${\rm
Comb}_P(\mathcal{A}_i)_{i\in I\cup\{\infty\}}$\index{${\rm
Comb}_P(\mathcal{A}_i)_{i\in I\cup\{\infty\}}$} for ${\rm
Comb}_P(\mathcal{A}_i)_{i\in I}$ with the empty structure
$\mathcal{A}_\infty$.

Note that if all predicates $P_i$ are disjoint, a structure
$\mathcal{A}_P$ is a $P$-combination and a disjoint union of
structures $\mathcal{A}_i$. In this case the $P$-combination
$\mathcal{A}_P$ is called {\em
disjoint}.\index{$P$-combination!disjoint} Clearly, for any
disjoint $P$-combination $\mathcal{A}_P$, ${\rm
Th}(\mathcal{A}_P)={\rm Th}(\mathcal{A}'_P)$, where
$\mathcal{A}'_P$ is obtained from $\mathcal{A}_P$ replacing
$\mathcal{A}_i$ by pairwise disjoint
$\mathcal{A}'_i\equiv\mathcal{A}_i$, $i\in I$. Thus, in this case,
similar to structures the $P$-operator works for the theories
$T_i={\rm Th}(\mathcal{A}_i)$ producing the theory $T_P={\rm
Th}(\mathcal{A}_P)$\index{$T_P$}, being {\em
$P$-combination}\index{$P$-combination} of $T_i$, which is denoted
by ${\rm Comb}_P(T_i)_{i\in I}$.\index{${\rm Comb}_P(T_i)_{i\in
I}$}

For an equivalence relation $E$ replacing disjoint predicates
$P_i$ by $E$-classes we get the structure
$\mathcal{A}_E$\index{$\mathcal{A}_E$} being the {\em
$E$-union}\index{$E$-union} of the structures $\mathcal{A}_i$. In
this case the operator mapping $(\mathcal{A}_i)_{i\in I}$ to
$\mathcal{A}_E$ is the {\em $E$-operator}\index{$E$-operator}. The
structure $\mathcal{A}_E$ is also called the {\em
$E$-combination}\index{$E$-combination} of the structures
$\mathcal{A}_i$ and denoted by ${\rm Comb}_E(\mathcal{A}_i)_{i\in
I}$\index{${\rm Comb}_E(\mathcal{A}_i)_{i\in I}$}; here
$\mathcal{A}_i=(\mathcal{A}_E\upharpoonright
A_i)\upharpoonright\Sigma(\mathcal{A}_i)$, $i\in I$. Similar
above, structures $\mathcal{A}'$, which are elementary equivalent
to $\mathcal{A}_E$, are denoted by ${\rm
Comb}_E(\mathcal{A}'_j)_{j\in J}$, where $\mathcal{A}'_j$ are
restrictions of $\mathcal{A}'$ to its $E$-classes. The
$E$-operator works for the theories $T_i={\rm Th}(\mathcal{A}_i)$
producing the theory $T_E={\rm Th}(\mathcal{A}_E)$\index{$T_E$},
being {\em $E$-combination}\index{$E$-combination} of $T_i$, which
is denoted by ${\rm Comb}_E(T_i)_{i\in I}$\index{${\rm
Comb}_E(T_i)_{i\in I}$} or by ${\rm
Comb}_E(\mathcal{T})$\index{${\rm Comb}_E(\mathcal{T})$}, where
$\mathcal{T}=\{T_i\mid i\in I\}$.

Clearly, $\mathcal{A}'\equiv\mathcal{A}_P$ realizing $p_\infty(x)$
is not elementary embeddable into $\mathcal{A}_P$ and can not be
represented as a disjoint $P$-combination of
$\mathcal{A}'_i\equiv\mathcal{A}_i$, $i\in I$. At the same time,
there are $E$-combinations such that all
$\mathcal{A}'\equiv\mathcal{A}_E$ can be represented as
$E$-combinations of some $\mathcal{A}'_j\equiv\mathcal{A}_i$. We
call this representability of $\mathcal{A}'$ to be the {\em
$E$-representability}.

If there is $\mathcal{A}'\equiv\mathcal{A}_E$ which is not
$E$-representable, we have the $E'$-representability replacing $E$
by $E'$ such that $E'$ is obtained from $E$ adding equivalence
classes with models for all theories $T$, where $T$ is a theory of
a restriction $\mathcal{B}$ of a structure
$\mathcal{A}'\equiv\mathcal{A}_E$ to some $E$-class and
$\mathcal{B}$ is not elementary equivalent to the structures
$\mathcal{A}_i$. The resulting structure $\mathcal{A}_{E'}$ (with
the $E'$-representability) is a {\em
$e$-completion}\index{$e$-completion}, or a {\em
$e$-saturation}\index{$e$-saturation}, of $\mathcal{A}_{E}$. The
structure $\mathcal{A}_{E'}$ itself is called {\em
$e$-complete}\index{Structure!$e$-complete}, or {\em
$e$-saturated}\index{Structure!$e$-saturated}, or {\em
$e$-universal}\index{Structure!$e$-universal}, or {\em
$e$-largest}\index{Structure!$e$-largest}.

For a structure $\mathcal{A}_E$ the number of {\em
new}\index{Structure!new} structures with respect to the
structures $\mathcal{A}_i$, i.~e., of the structures $\mathcal{B}$
which are pairwise elementary non-equivalent and elementary
non-equivalent to the structures $\mathcal{A}_i$, is called the
{\em $e$-spectrum}\index{$e$-spectrum} of $\mathcal{A}_E$ and
denoted by $e$-${\rm Sp}(\mathcal{A}_E)$.\index{$e$-${\rm
Sp}(\mathcal{A}_E)$} The value ${\rm sup}\{e$-${\rm
Sp}(\mathcal{A}'))\mid\mathcal{A}'\equiv\mathcal{A}_E\}$ is called
the {\em $e$-spectrum}\index{$e$-spectrum} of the theory ${\rm
Th}(\mathcal{A}_E)$ and denoted by $e$-${\rm Sp}({\rm
Th}(\mathcal{A}_E))$.\index{$e$-${\rm Sp}({\rm
Th}(\mathcal{A}_E))$}

If $\mathcal{A}_E$ does not have $E$-classes $\mathcal{A}_i$,
which can be removed, with all $E$-classes
$\mathcal{A}_j\equiv\mathcal{A}_i$, preserving the theory ${\rm
Th}(\mathcal{A}_E)$, then $\mathcal{A}_E$ is called {\em
$e$-prime}\index{Structure!$e$-prime}, or {\em
$e$-minimal}\index{Structure!$e$-minimal}.

For a structure $\mathcal{A}'\equiv\mathcal{A}_E$ we denote by
${\rm TH}(\mathcal{A}')$ the set of all theories ${\rm
Th}(\mathcal{A}_i)$\index{${\rm Th}(\mathcal{A}_i)$} of
$E$-classes $\mathcal{A}_i$ in $\mathcal{A}'$.

By the definition, an $e$-minimal structure $\mathcal{A}'$
consists of $E$-classes with a minimal set ${\rm
TH}(\mathcal{A}')$. If ${\rm TH}(\mathcal{A}')$ is the least for
models of ${\rm Th}(\mathcal{A}')$ then $\mathcal{A}'$ is called
{\em $e$-least}.\index{Structure!$e$-least}

\medskip
{\bf Definition} \cite{cl}. Let $\overline{\mathcal{T}}$ be the
class of all complete elementary theories of relational languages.
For a set $\mathcal{T}\subset\overline{\mathcal{T}}$ we denote by
${\rm Cl}_E(\mathcal{T})$ the set of all theories ${\rm
Th}(\mathcal{A})$, where $\mathcal{A}$ is a structure of some
$E$-class in $\mathcal{A}'\equiv\mathcal{A}_E$,
$\mathcal{A}_E={\rm Comb}_E(\mathcal{A}_i)_{i\in I}$, ${\rm
Th}(\mathcal{A}_i)\in\mathcal{T}$. As usual, if $\mathcal{T}={\rm
Cl}_E(\mathcal{T})$ then $\mathcal{T}$ is said to be {\em
$E$-closed}.\index{Set!$E$-closed}

The operator ${\rm Cl}_E$ of $E$-closure can be naturally extended
to the classes $\mathcal{T}\subset\overline{\mathcal{T}}$ as
follows: ${\rm Cl}_E(\mathcal{T})$ is the union of all ${\rm
Cl}_E(\mathcal{T}_0)$ for subsets
$\mathcal{T}_0\subseteq\mathcal{T}$.

For a set $\mathcal{T}\subset\overline{\mathcal{T}}$ of theories
in a language $\Sigma$ and for a sentence $\varphi$ with
$\Sigma(\varphi)\subseteq\Sigma$ we denote by
$\mathcal{T}_\varphi$\index{$\mathcal{T}_\varphi$} the set
$\{T\in\mathcal{T}\mid\varphi\in T\}$.

\medskip
{\bf Proposition 1.1} \cite{cl}. {\em If
$\mathcal{T}\subset\overline{\mathcal{T}}$ is an infinite set and
$T\in\overline{\mathcal{T}}\setminus\mathcal{T}$ then $T\in{\rm
Cl}_E(\mathcal{T})$ {\rm (}i.e., $T$ is an {\sl accumulation
point} for $\mathcal{T}$ with respect to $E$-closure ${\rm
Cl}_E${\rm )} if and only if for any formula $\varphi\in T$ the
set $\mathcal{T}_\varphi$ is infinite.}

\medskip
{\bf Theorem 1.2} \cite{cl}. {\em For any sets
$\mathcal{T}_0,\mathcal{T}_1\subset\overline{\mathcal{T}}$, ${\rm
Cl}_E(\mathcal{T}_0\cup\mathcal{T}_1)={\rm
Cl}_E(\mathcal{T}_0)\cup{\rm Cl}_E(\mathcal{T}_1)$.}

\medskip
{\bf Definition} \cite{cl}. Let $\mathcal{T}_0$ be a closed set in
a topological space $(\mathcal{T},\mathcal{O}_E(\mathcal{T}))$,
where $\mathcal{O}_E(\mathcal{T})=\{\mathcal{T}\setminus{\rm
Cl}_E(\mathcal{T}')\mid\mathcal{T}'\subseteq\mathcal{T}\}$. A
subset $\mathcal{T}'_0\subseteq\mathcal{T}_0$ is said to be {\em
generating}\index{Set!generating} if $\mathcal{T}_0={\rm
Cl}_E(\mathcal{T}'_0)$. The generating set $\mathcal{T}'_0$ (for
$\mathcal{T}_0$) is {\em minimal}\index{Set!generating!minimal} if
$\mathcal{T}'_0$ does not contain proper generating subsets. A
minimal generating set $\mathcal{T}'_0$ is {\em
least}\index{Set!generating!least} if $\mathcal{T}'_0$ is
contained in each generating set for $\mathcal{T}_0$.

\medskip
{\bf Theorem 1.3} \cite{cl}. {\em If $\mathcal{T}'_0$ is a
generating set for a $E$-closed set $\mathcal{T}_0$ then the
following conditions are equivalent:

$(1)$ $\mathcal{T}'_0$ is the least generating set for
$\mathcal{T}_0$;

$(2)$ $\mathcal{T}'_0$ is a minimal generating set for
$\mathcal{T}_0$;

$(3)$ any theory in $\mathcal{T}'_0$ is isolated by some set
$(\mathcal{T}'_0)_\varphi$, i.e., for any $T\in\mathcal{T}'_0$
there is $\varphi\in T$ such that
$(\mathcal{T}'_0)_\varphi=\{T\}$;

$(4)$ any theory in $\mathcal{T}'_0$ is isolated by some set
$(\mathcal{T}_0)_\varphi$, i.e., for any $T\in\mathcal{T}'_0$
there is $\varphi\in T$ such that
$(\mathcal{T}_0)_\varphi=\{T\}$.}

\section{Relative $e$-spectra and their properties}

{\bf Definition.} For a structure $\mathcal{A}_E$ and a class $K$
of structures, the number of {\em new}\index{Structure!new}
structures with respect to the structures $\mathcal{A}_i$ and to
the class $K$, i.~e., of the structures $\mathcal{B}$ forming
$E$-classes of models of ${\rm Th}(\mathcal{A}_E)$ such that
$\mathcal{B}$ are pairwise elementary non-equivalent and
elementary non-equivalent to the structures $\mathcal{A}_i$ in
$\mathcal{A}_E$ as well as to the structures in $K$, is called the
{\em relative $e$-spectrum}\index{$e$-spectrum!relative} of
$\mathcal{A}_E$ with respect to $K$ and denoted by $e_K$-${\rm
Sp}(\mathcal{A}_E)$.\index{$e_K$-${\rm Sp}(\mathcal{A}_E)$} The
value ${\rm sup}\{e_K$-${\rm
Sp}(\mathcal{A}'))\mid\mathcal{A}'\equiv\mathcal{A}_E\}$ is called
the {\em relative $e$-spectrum}\index{$e$-spectrum!relative} of
the theory ${\rm Th}(\mathcal{A}_E)$ with respect to $K$ and
denoted by $e_K$-${\rm Sp}({\rm
Th}(\mathcal{A}_E))$.\index{$e_K$-${\rm Sp}({\rm
Th}(\mathcal{A}_E))$}

Similarly for a class $\mathcal{T}$ of theories and for a theory
$T={\rm Th}(\mathcal{A}_E)$ we denote by $e_\mathcal{T}$-${\rm
Sp}(T)$\index{$e_\mathcal{T}$-${\rm Sp}(T)$} the value $e_K$-${\rm
Sp}(T)$, where $K=K(\mathcal{T})$ is the class of all structures,
each of which is a model of a theory in $\mathcal{T}$. The value
$e_\mathcal{T}$-${\rm Sp}(T)$ is called the {\em relative
$e$-spectrum}\index{$e$-spectrum!relative} of the theory $T$ with
respect to $\mathcal{T}$.

\medskip
{\bf Remark 2.1.} 1. the class $K(\mathcal{T})$, in the definition
above, can be replaced by any subclass $K'\subseteq
K(\mathcal{T})$ such that any structure in $K(\mathcal{T})$ is
elementary equivalent to a structure in $K'$.

2. if $K_1\subseteq K_2$ then $e_{K_1}$-${\rm Sp}(T)\geq
e_{K_2}$-${\rm Sp}(T)$, and if $\mathcal{T}_1\subseteq
\mathcal{T}_2$ then $e_{\mathcal{T}_1}$-${\rm Sp}(T)\geq
e_{\mathcal{T}_2}$-${\rm Sp}(T)$.

3. The value $e_\mathcal{T}$-${\rm Sp}(T)$ is equal to the
supremum $|\mathcal{T}_1\setminus\mathcal{T}_0|$ for theories of
$E$-classes of models of $T$ such that $\mathcal{T}_1$ consists of
all these theories and $\mathcal{T}_0\subseteq\mathcal{T}_1$ with
${\rm Cl}_E(\mathcal{T}_0)=\mathcal{T}_1$.

\medskip
{\bf Definition.} Two theories $T_1$ and $T_2$ of a language
$\Sigma$ are {\em disjoint}\index{Theories!disjoint modulo
$\Sigma_0$} modulo $\Sigma_0$, where $\Sigma_0\subseteq\Sigma$, or
{\em $\Sigma_0$-disjoint}\index{Theories!$\Sigma_0$-disjoint} if
$T_1$ and $T_2$ are do not have common nonempty predicates for
$\Sigma\setminus\Sigma_0$. If $T_1$ and $T_2$ are
$\varnothing$-disjoint, these theories are called simply {\em
disjoint}\index{Theories!disjoint}.

Families $\mathcal{T}_j$, $j\in J$, of theories in the language
$\Sigma$ are {\em disjoint}\index{Families!disjoint} modulo
$\Sigma_0$, or {\em
$\Sigma_0$-disjoint}\index{Families!$\Sigma_0$-disjoint} if
$T_{j_1}$ and $T_{j_2}$ are $\Sigma_0$-disjoint for any
$T_{j_1}\in\mathcal{T}_{j_1}$, $T_{j_2}\in\mathcal{T}_{j_2}$,
$j_1\ne j_2$. If $T_{j_1}$ and $T_{j_2}$ are disjoint for any
$T_{j_1}\in\mathcal{T}_{j_1}$, $T_{j_2}\in\mathcal{T}_{j_2}$,
$j_1\ne j_2$, then the families $\mathcal{T}_j$, $j\in J$, are
{\em disjoint}\index{Families!disjoint} too.

\medskip
The following properties are obvious.

1. Any families of theories in a language $\Sigma$ are
$\Sigma$-disjoint.

2. (Monotony) If $\Sigma_0\subseteq\Sigma_1\subseteq\Sigma$ then
disjoint families modulo $\Sigma_0$, in the language $\Sigma$, are
disjoint modulo $\Sigma_1$.

3. (Monotony) If families $\mathcal{T}_{j_1}$ and
$\mathcal{T}_{j_2}$ are $\Sigma_0$-disjoint then any subfamilies
$\mathcal{T}'_{j_1}\subseteq\mathcal{T}_{j_1}$ and
$\mathcal{T}'_{j_2}\subseteq\mathcal{T}_{j_2}$ are
$\Sigma_0$-disjoint too.

\medskip
Below we denote by $K_\Sigma$\index{$K_\Sigma$} the class of all
structures in languages containing $\Sigma$ such that all
predicates outside $\Sigma$ are empty. Similarly we denote by
$\mathcal{T}_\Sigma$\index{$\mathcal{T}_\Sigma$} the class of all
theories of structures in $K_\Sigma$.

\medskip
{\bf Theorem 2.2.} (Relative additivity for $e$-spectra) {\em If
$\mathcal{T}_j$, $j\in J$, are $\Sigma_0$-disjoint
families then for the $E$-combination
$T={\rm Comb}_E(T_i)_{i\in I}$ of $\{T_i\mid i\in
I\}=\bigcup\limits_{j\in J}\mathcal{T}_j$ and for the
$E$-combinations $T_j={\rm Comb}_E(\mathcal{T}_j)$, $j\in J$,
\begin{equation}\label{rest1}
e_{\mathcal{T}_{\Sigma_0}}\mbox{-}{\rm Sp}(T)=\sum\limits_{j\in
J}(e_{\mathcal{T}_{\Sigma_0}}\mbox{-}{\rm
Sp}(T_j)).\end{equation}}

\medskip
{\bf\em Proof.} Denote by $\mathcal{T}$ the set of theories for
$E$-classes of models of $T$. Since the families $\mathcal{T}_j$
are $\Sigma_0$-disjoint, 
applying
Proposition 1.1 we have that a theory $T^\ast$ belongs to ${\rm
Cl}_E(\mathcal{T}^\ast)$, where
$\mathcal{T}^\ast\subseteq\mathcal{T}$, if and only if some of the
following conditions holds:

1) $T^\ast\in\mathcal{T}^\ast$;

2) for any formula $\varphi\in T^\ast$ without predicate symbols
in $\Sigma\setminus\Sigma_0$, or with predicate symbols in
$\Sigma\setminus\Sigma_0$ and saying that corresponding predicates
are empty, there are infinitely many theories in
$T\in\mathcal{T}^\ast$ containing $\varphi$;

3) for any formula $\varphi\in T^\ast$, saying that some
predicates in $\Sigma\setminus\Sigma_0$ which used in $\varphi$
are nonempty, there are infinitely many theories in
$T\in\mathcal{T}^\ast\cap\mathcal{T}_j$, for some $j$, containing
$\varphi$; moreover, the theories $T$ belong to the unique
$\mathcal{T}_j$.

Indeed, taking a formula $\varphi$ in the language $\Sigma$ we
have finitely many symbols $R_1,\ldots,R_n$ in
$\Sigma\setminus\Sigma_0$, used in $\varphi$. Considering formulas
$\psi_i$ saying that $R_k$ are nonempty, $k=1,\ldots,n$, we get
finitely many possibilities for
$\chi^{\delta_1,\ldots,\delta_n}\rightleftharpoons\varphi\wedge\bigwedge\limits_{k=1}^n\psi^{\delta_k}_k$,
$\delta_k\in\{0,1\}$. Since $\varphi$ is equivalent to
$\bigvee\limits_{\delta_1,\ldots,\delta_n}\chi^{\delta_1,\ldots,\delta_n}$
and only subdisjunctions with positive $\psi_k$ related to the
fixed $\mathcal{T}_j$ hold, we can divide the disjunction to
disjoint parts related to $\mathcal{T}_j$. Since for $\varphi$
there are finitely many related $\mathcal{T}_j$, we have finitely
many cases for $\varphi$, each of which related to the fixed
$\mathcal{T}_j$. These cases are described in Item 3. Item 2 deals
with formulas in the language $\Sigma_0$ and with formulas for
empty part in $\Sigma\setminus\Sigma_0$. In particular, by
Proposition 1.1 these formulas define ${\rm
Cl}_E(\mathcal{T}^\ast)\cap\mathcal{T}_{\Sigma_0}$.

Using Items 1--3 we have for $\mathcal{T}^\ast$ that a theory
$T^\ast$ belongs to
$\mathcal{T}^\ast\setminus\mathcal{T}_{\Sigma_0}$ if and only if
$T^\ast$ belong to
$(\mathcal{T}^\ast\cap\mathcal{T}_j)\setminus\mathcal{T}_{\Sigma_0}$
for unique $j\in J$. Thus theories witnessing the value
$e_{\mathcal{T}_{\Sigma_0}}\mbox{-}{\rm Sp}(T)$ are divided into
disjoint parts witnessing the values
$e_{\mathcal{T}_{\Sigma_0}}\mbox{-}{\rm Sp}(T_j)$. Thus the
equality (\ref{rest1}) holds.~$\Box$

\medskip
{\bf Remark 2.3.} Having positive ${\rm ComLim}$ \cite{cs} the
equality (\ref{rest1}) can fail if families $\mathcal{T}_j$ are
not $\Sigma_0$-disjoint, even for finite sets $J$ of indexes,
producing
\begin{equation}\label{rest2}
e_{\mathcal{T}_{\Sigma_0}}\mbox{-}{\rm Sp}(T')<\sum\limits_{j\in
J}(e_{\mathcal{T}_{\Sigma_0}}\mbox{-}{\rm Sp}(T_j))
\end{equation} for
appropriate $T'$.

\medskip
Theorem 2.2 immediately implies

\medskip
{\bf Corollary 2.4.} {\em If $\mathcal{T}_j$, $j\in J$, are
disjoint then for the $E$-combination $T={\rm Comb}_E(T_i)_{i\in
I}$ of $\{T_i\mid i\in I\}=\bigcup\limits_{j\in J}\mathcal{T}_j$
and for the $E$-combinations $T_j={\rm Comb}_E(\mathcal{T}_j)$,
$j\in J$,
\begin{equation}\label{rest3}
e_{\mathcal{T}_\varnothing}\mbox{-}{\rm Sp}(T)=\sum\limits_{j\in
J}(e_{\mathcal{T}_\varnothing}\mbox{-}{\rm
Sp}(T_j)).\end{equation}}

\medskip
{\bf Definition.} The theory $T$ in Theorem 2.2 is called the {\em
$\Sigma_0$-disjoint 
$E$-union}\index{$E$-union of theories!$\Sigma_0$-disjoint}
of the theories $T_j$, $j\in J$, and the
theory $T$ in Corollary 2.4 is the {\em disjoint
$E$-union}\index{$E$-union of theories!disjoint} of the theories
$T_j$, $j\in J$.

\medskip
{\bf Remark 2.5.} Additivity (\ref{rest1}) and, in particular,
(\ref{rest3}) can be failed without indexes
$\mathcal{T}_{\Sigma_0}$. Indeed, it is possible to find $T_j$
with $e\mbox{-}{\rm Sp}(T_j)=0$ (for instance, with finite
$\mathcal{T}_j$) while $e\mbox{-}{\rm Sp}(T)$ can be positive.
Take, for example, disjoint singletons $\mathcal{T}_n=\{T_n\}$,
$n\in\omega\setminus\{0\}$, such that $T_n$ has $n$-element
models. We have $e\mbox{-}{\rm Sp}(T_n)=0$ for each $n$ while
$e\mbox{-}{\rm Sp}(T)=1$, since the theory
$T_\infty\in\mathcal{T}_\varnothing$ with infinite models belong
to ${\rm Cl}_E(\{T_n\mid n\in\omega\setminus\{0\}\})$. Thus, for
disjoint families $\mathcal{T}_j$, $j\in J$, the equality
\begin{equation}\label{rest4}
e\mbox{-}{\rm Sp}(T)=\sum\limits_{j\in J}(e\mbox{-}{\rm
Sp}(T_j))\end{equation} can fail. Moreover, producing the effect
above for definable subsets in models of $T_j$ we get
$$e_{\mathcal{T}_{\Sigma_0}}\mbox{-}{\rm Sp}(T)>\sum\limits_{j\in
J}(e_{\mathcal{T}_{\Sigma_0}}\mbox{-}{\rm Sp}(T_j)).$$

At the same time, by Corollary 2.4 (respectively, by Theorem 2.2)
the equality (\ref{rest4}) holds for ($\Sigma_0$-)disjoint
families $\mathcal{T}_j$, $j\in J$, if $J$ is finite and each
$\mathcal{T}_j$ does not generate theories in
$\mathcal{T}_\varnothing$ (in $\mathcal{T}_{\Sigma_0}$).

\medskip
Applying the equality (\ref{rest3}) we take an $E$-combination
$T_0$ with $e_{\mathcal{T}_\varnothing}\mbox{-}{\rm
Sp}(T_0)=\lambda$. Furthermore we consider disjoint copies $T_j$,
$j\in J$, of $T_0$. Combining $E$-classes of all $T_j$ we obtain a
theory $T$ such that if $J$ is finite then
$e_{\mathcal{T}_\varnothing}\mbox{-}{\rm Sp}(T)=|J|\cdot\lambda$.
We have the same formula if $|J|\geq\omega$ and $\lambda>0$ since,
in this case, the $E$-closure for theories of $E$-classes of
models of $T$ consists of theories of $E$-classes for theories
$T_j$ as well some theories in $\mathcal{T}_\varnothing$. If
$E$-classes have a fixed finite or only infinite cardinalities,
this theory has models whose cardinalities (finite or countable)
are equal to the (either finite or countable) cardinality of
models of $T_j$. Similarly, having theories $T_\lambda$ of
languages $\Sigma$ with cardinalities $|\Sigma|=\lambda+1$ and
with $e\mbox{-}{\rm Sp}(T_0)=\lambda>0$ \cite[Proposition 4.3]{cs}
and taking $E$-combinations with their disjoint copies we get

\medskip
{\bf Proposition 2.6.} {\em For any positive cardinality $\lambda$
there is a theory $T$ such that $E$-classes of models of $T$ form
copies $T_j$, $j\in J$, of some $E$-combination $T_0$ with a
language $\Sigma$ in the cardinality $\lambda+1$, with
$e_{\mathcal{T}_\varnothing}\mbox{-}{\rm Sp}(T_0)=\lambda$, and
$e_{\mathcal{T}_\varnothing}\mbox{-}{\rm Sp}(T)=|J|\cdot\lambda$.}

\medskip
{\bf Remark 2.7.} Since there are required theories $T_0$ which do
not generate $E$-classes for $\mathcal{T}_\varnothing$,
Proposition 2.6 can be reformulated without the index
$\mathcal{T}_\varnothing$.

\medskip
{\bf Remark 2.8.} Extending the $\Sigma_0$-disjoint
$\Sigma_0$-coordinated $E$-union $T$ by definable bijections
linking $E$-classes we can omit the additivity (\ref{rest1}).
Indeed, adding, for instance, bijections $f_{jk}$ witnessing
isomorphisms for models of disjoint copies $T_j$ and $T_j$, have
we $e_{\mathcal{T}_\varnothing}\mbox{-}{\rm Sp}(T_j)$ instead of
$e_{\mathcal{T}_\varnothing}\mbox{-}{\rm
Sp}(T_j)+e_{\mathcal{T}_\varnothing}\mbox{-}{\rm Sp}(T_k)$. Thus,
bijections $f_{jk}$ allow to vary
$e_{\mathcal{T}_\varnothing}\mbox{-}{\rm Sp}(T)$ from $\lambda$ to
$|J|\cdot\lambda$ in terms of Proposition 2.6. Thus the equality
(\ref{rest1}) can fail again producing (\ref{rest2}) for
appropriate~$T'$.

\section{Families of theories with(out) least generating sets}

Below we apply Theorem 1.3 characterizing the existence of
$e$-least generating sets for $\Sigma_0$-disjoint families of
theories.

The following natural questions arises:

\medskip
{\bf Question 1.} {\em When the existence of the least generating
sets for the families $\mathcal{T}_j$, $j\in J$, is equivalent to
the existence of the least generating set for the family
$\bigcup\limits_{j\in J}\mathcal{T}_j$?}

\medskip
{\bf Question 2.} {\em Is it true that under conditions of Theorem
$2.2$ the existence of the least generating sets for the families
$\mathcal{T}_j$, $j\in J$, is equivalent to the existence of the
least generating set for the family $\bigcup\limits_{j\in
J}\mathcal{T}_j$?}

\medskip
Considering Question 2, we note below that the property of the
(non)existence of the least generating sets is not preserving
under expansions and extensions of families of theories.

\medskip
{\bf Proposition 3.1.} {\em Any $E$-closed family $\mathcal{T}_0$
of theories in a language $\Sigma_0$ can be transformed to an
$E$-closed family $\mathcal{T}'_0$ in a language
$\Sigma'_0\supseteq\Sigma_0$ such that $\mathcal{T}'_0$ consists
of expansions of theories in $\mathcal{T}_0$ and $\mathcal{T}'_0$
has the least generating set.}

\medskip
{\bf\em Proof.} Forming $\Sigma'_0$ it suffices to take new
predicate symbols $R_{T_0}$, $T_0\in\mathcal{T}_0$, such that
$R_{T_0}\ne\varnothing$ for interpretations in the models of
expansion $T'_0$ of $T_0$ and $R_{T_0}=\varnothing$ for
interpretations in the models of expansion $T'_1$ of $T_1\ne T_0$.
Each formula $\exists\bar{x}R_{T_0}(\bar{x})$ isolates $T'_0$, and
thus $\mathcal{T}'_0$ has the least generating set in view of
Theorem 1.3.~$\Box$

\medskip
Existence of families $\mathcal{T}_0$ without least generating
sets implies

\medskip
{\bf Corollary 3.2.} {\em The property of non-existence of least
generating sets is not preserved under expansions of theories.}

\medskip
{\bf Remark 3.3.} The expansion $\mathcal{T}'_0$ of
$\mathcal{T}_0$ in the proof of Proposition 3.1 produces discrete
topologies for sets of theories in
$\mathcal{T}_0\cup\mathcal{T}'_0$. In fact, for this purpose it
suffices to isolate finite sets in $\mathcal{T}_0$ since any two
distinct elements $T_0,T_1\in\mathcal{T}_0$ are separated by
formulas $\varphi$ such that $\varphi\in T_i$ and $\neg\varphi\in
T_{1-i}$, $i=0,1$.

Note also that these operators of discretization transform the
given set $\mathcal{T}_0$ to a set $\mathcal{T}'_0$ with identical
${\rm Cl}_E$.

\medskip
Clearly, if a set $\mathcal{T}_0$ has the discrete topology it can
not be expanded to a set without the least generating set. At the
same time, there are expansions that transform sets with the least
generating sets to sets without the least generating sets. Indeed,
take Example in \cite[Remark 3]{lut} with countably many disjoint
copies $\mathcal{F}_q$, $q\in\mathbb Q$, of linearly ordered sets
isomorphic to $\langle\omega,\leq\rangle$ and ordering limits
$J_q=\overline{\rm lim}\,F_q$ by the ordinary dense order on
$\mathbb Q$ such that $\{J_q\mid q\in\mathbb Q\}$ is densely
ordered. We have a dense interval $\{J_q\mid q\in\mathbb Q\}$
whereas the set $\cup\{F_q\mid q\in\mathbb Q\}$ forms the least
generating set $\mathcal{T}_0$ of theories for ${\rm
Cl}_E(\mathcal{T}_0)$. Now we expand the ${\rm LU}$-theories for
$\mathcal{F}_q$ and $J_q$ by new predicate symbol $R$ such that
$R$ is empty for all theories corresponding to $\mathcal{F}_q$ and
$\forall\bar{x}R(\bar{x})$ is satisfied for all theories
corresponding to $J_q$. The predicate $R$ separates the set of
theories for $J_q$ with respect to ${\rm Cl}_E$. At the same time
the theories for $J_q$ forms the dense interval producing the set
without the least generating set in view of \cite[Theorem 2]{lut}.
Thus, we get the following

\medskip
{\bf Proposition 3.4.} {\em There is an $E$-closed family
$\mathcal{T}_0$ of theories in a language $\Sigma_0$ and with the
least generating set, which can be transformed to an $E$-closed
family $\mathcal{T}'_0$ in a language $\Sigma'_0\supseteq\Sigma_0$
such that $\mathcal{T}'_0$ consists of expansions of theories in
$\mathcal{T}_0$ and $\mathcal{T}'_0$ does not have the least
generating set.}

\medskip
{\bf Corollary 3.5.} {\em The property of existence of least
generating sets is not preserved under expansions of theories.}

\medskip
{\bf Remark 3.6.} Adding the predicate $R$ which separates
theories for $J_q$ from theories for $\mathcal{F}_q$, we get a
copy for each $J_q$ containing empty $R$. This effect is based on
the property that separating an accumulation point $J_q$ for
$\mathcal{F}_q$ we get new accumulation point preserving formulas
in the initial language.

Introducing the predicate $R$ together with the discretization for
$\mathcal{F}_q$, $E$-closures do not generate new theories.

\medskip
{\bf Proposition 3.7.} {\em Any family $\mathcal{T}_0$ of theories
in a language $\Sigma$, with infinitely many empty predicates for
all theories in $\mathcal{T}_0$, can be extended to a family
$\mathcal{T}'_0$ in the language $\Sigma$ such that
$\mathcal{T}'_0$ does not have the least generating set.}

\medskip
{\bf\em Proof.} Let $\Sigma_0\subseteq\Sigma$ consist of predicate
symbols which are empty for all theories in $\mathcal{T}_0$. Now
we consider a family $\mathcal{T}_1$ of ${\rm LU}$-theories such
that all these theories have empty predicates for
$\Sigma\setminus\Sigma_0$, and, using $\Sigma_0$ as for
\cite[Theorem 2]{lut}, $\mathcal{T}_1$ does not have the least
generating set forming a dense interval. The family
$\mathcal{T}'_0=\mathcal{T}_0\,\dot{\cup}\,\mathcal{T}_1$ extends
$\mathcal{T}_0$ and does not have the least generating set since
for any $\mathcal{T}''_0\subseteq\mathcal{T}'_0$, ${\rm
Cl}_E(\mathcal{T}''_0)={\rm
Cl}_E(\mathcal{T}''_0\cap\mathcal{T}_0)\,\dot{\cup}\,{\rm
Cl}_E(\mathcal{T}''_0\cap\mathcal{T}_1)$.~$\Box$

\medskip
{\bf Corollary 3.8.} {\em The property of existence of least
generating sets is not preserved under extensions of sets of
theories.}

\medskip
In view of Theorem 1.3 any family consisting of all theories in a
given infinite language both does not have the least generating
set and does not have a proper extension in the given language.
Thus there are families of theories without least generating sets
and without extension having least generating sets. At the same
time the following proposition holds.

\medskip
{\bf Proposition 3.9.} {\em There is an $E$-closed family
$\mathcal{T}_0$ of theories in a language $\Sigma$ and without the
least generating set such that $\mathcal{T}_0$ can be extended to
an $E$-closed family $\mathcal{T}'_0$ in the language $\Sigma$ and
with the least generating set.}

\medskip
{\bf\em Proof.} It suffices to take Example in \cite[Remark
3]{lut} that we used for the proof of Proposition 3.4. The
theories for $\{J_q\mid q\in\mathbb Q\}$ form a family without the
least generating set whereas an extension of this family by the
theories for $\mathcal{F}_q$ has the least generating set.~$\Box$

\medskip
{\bf Corollary 3.10.} {\em The property of non-existence of least
generating sets is not preserved under extensions of sets of
theories.}

\medskip
{\bf Remark 3.11.} If an extension of an $E$-closed family
$\mathcal{T}_0$ of theories transforms $\mathcal{T}_0$ with the
least generating set to an $E$-closed family $\mathcal{T}'_0$
without the least generating set then, in view of Theorem 1.3,
having the generating set in $\mathcal{T}_0$ consisting of
isolated points we lose this property for $\mathcal{T}'_0$. If an
extension of an $E$-closed family $\mathcal{T}_0$ of theories
transforms $\mathcal{T}_0$ without the least generating set to an
$E$-closed family $\mathcal{T}'_0$ with the least generating set
then, again in view of Theorem 1.3, we add a set of isolated
theories to $\mathcal{T}_0$ generating all theories in
$\mathcal{T}'_0$.

\medskip
Now we return to Questions 1 and 2.

Clearly, for any set $\mathcal{T}$ of theories, ${\rm
Cl}_E(\mathcal{T}\cap\mathcal{T}_{\Sigma_0})\subset\mathcal{T}_{\Sigma_0}$.
Therefore ${\rm Cl}_E(\mathcal{T})$ and each its generating set
are divided into parts: in $\mathcal{T}_{\Sigma_0}$ and disjoint
with $\mathcal{T}_{\Sigma_0}$. Since $\mathcal{T}_j$, $j\in J$,
are disjoint with respect to $\mathcal{T}_{\Sigma_0}$, each
$\mathcal{T}_j$ has the least generating set if and only if both
$\mathcal{T}_j\cap\mathcal{T}_{\Sigma_0}$ and
$\mathcal{T}_j\setminus\mathcal{T}_{\Sigma_0}$ have the least
generating sets. Since under conditions of Theorem 2.2 the sets
$\mathcal{T}_j\setminus\mathcal{T}_{\Sigma_0}$ are disjoint, $j\in
J$, we have the following proposition answering Question 1.

\medskip
{\bf Proposition 3.12.} {\em The set $\bigcup\limits_{j\in
J}\mathcal{T}_j$ has the least generating set if and only if
$\left(\bigcup\limits_{j\in
J}\mathcal{T}_j\right)\cap\mathcal{T}_{\Sigma_0}$ has the least
generating set and each
$\mathcal{T}_j\setminus\mathcal{T}_{\Sigma_0}$ has the least
generating set.}

\medskip
Since $\left(\bigcup\limits_{j\in
J}\mathcal{T}_j\right)\cap\mathcal{T}_{\Sigma_0}$ can be an
arbitrary extension of each
$\mathcal{T}_j\cap\mathcal{T}_{\Sigma_0}$, Propositions 3.7 and
3.12 imply the following corollary answering Question~2.

\medskip
{\bf Corollary 3.13.} {\em For any infinite language $\Sigma_0$
there are $\Sigma_0$-disjoint families $\mathcal{T}_j$, $j\in J$,
with the least generating sets such that $\bigcup\limits_{j\in
J}\mathcal{T}_j$ does not have the least generating set.}

\section{Relative closures and relative least generating sets}

{\bf Definition.} Let $\mathcal{T}$ be a class of theories. For a
set $\mathcal{T}_0\subset\overline{\mathcal{T}}$ we denote by
${\rm Cl}_{E,\mathcal{T}}(\mathcal{T}_0)$\index{${\rm
Cl}_{E,\mathcal{T}}(\mathcal{T}_0)$} the set ${\rm
Cl}_{E}(\mathcal{T}_0)\setminus\mathcal{T}$. The set ${\rm
Cl}_{E,\mathcal{T}}(\mathcal{T}_0)$ is called the {\em relative
$E$-closure}\index{$E$-closure!relative} of the set
$\mathcal{T}_0$ with respect to $\mathcal{T}$, or {\em
$\mathcal{T}$-relative
$E$-closure}.\index{$E$-closure!$\mathcal{T}$-relative} If
$\mathcal{T}_0\setminus\mathcal{T}={\rm
Cl}_{E,\mathcal{T}}(\mathcal{T}_0)$ then $\mathcal{T}_0$ is said
to be (relatively) {\em $E$-closed with respect to
$K$},\index{Set!$E$-closed!with respect to
$K$}\index{Set!relatively $E$-closed} or {\em
$\mathcal{T}$-relatively
$E$-closed}.\index{Set!$\mathcal{T}$-relatively $E$-closed}

Let $\mathcal{T}_0$ be a closed set in a topological space
$(\mathcal{T},\mathcal{O}_{E}(\mathcal{T}))$. A subset
$\mathcal{T}'_0\subseteq\mathcal{T}_0$ is said to be {\em
generating}\index{Set!generating!with respect to $\mathcal{T}$}
with respect to $\mathcal{T}$, or {\em $\mathcal{T}$-relatively
generating},\index{Set!generating!$\mathcal{T}$-relatively} if
$\mathcal{T}_0\setminus\mathcal{T}={\rm
Cl}_{E,\mathcal{T}}(\mathcal{T}'_0)$. The $\mathcal{T}$-relatively
generating set $\mathcal{T}'_0$ (for $\mathcal{T}_0$) is {\em
$\mathcal{T}$-minimal}\index{Set!$\mathcal{T}$-relatively
generating!$\mathcal{T}$-minimal} if
$\mathcal{T}'_0\setminus\mathcal{T}$ does not contain proper
subsets $\mathcal{T}''_0$ such that
$\mathcal{T}_0\setminus\mathcal{T}={\rm
Cl}_{E,\mathcal{T}}((\mathcal{T}'_0\cap\mathcal{T})\cup\mathcal{T}''_0)$.
A $\mathcal{T}$-minimal $\mathcal{T}$-relatively generating set
$\mathcal{T}'_0$ is {\em
$\mathcal{T}$-least}\index{Set!generating!$\mathcal{T}$-least} if
$\mathcal{T}'_0\setminus\mathcal{T}$ is contained in
$\mathcal{T}''_0\setminus\mathcal{T}$ for each
$\mathcal{T}$-relatively generating set $\mathcal{T}''_0$
for~$\mathcal{T}_0$.

\medskip
{\bf Remark 4.1.} Note that for $\mathcal{T}$-least generating
sets $\mathcal{T}'_0$, in general, we can say that
$\mathcal{T}'_0$ are uniquely defined only with respect to
$\mathcal{T}$. Moreover, since ${\rm
Cl}_E(\mathcal{T}_0\cup\mathcal{T}_1)={\rm
Cl}_E(\mathcal{T}_0)\cup{\rm Cl}_E(\mathcal{T}_1)$ for any sets
$\mathcal{T}_0,\mathcal{T}_1\subset\overline{\mathcal{T}}$ by
Theorem 1.2, then for $E$-closed $\mathcal{T}$, ${\rm
Cl}_E(\mathcal{T}'_0\cup\mathcal{T})={\rm
Cl}_E(\mathcal{T}'_0)\cup \mathcal{T}$ and $\mathcal{T}'_0$ is a
$\mathcal{T}$-least generating set if and only if
$\mathcal{T}'_0\cup \mathcal{T}'$ is a $\mathcal{T}$-least
generating set for some (any) $\mathcal{T}'\subseteq\mathcal{T}$,
as well as if and only if  $\mathcal{T}'_0\setminus\mathcal{T}$ is
a $\mathcal{T}$-least generating set.

\medskip
The following theorem generalizes Theorem 1.3.

\medskip
{\bf Theorem 4.2.} {\em If $\mathcal{T}$ is a $E$-closed set and
$\mathcal{T}'_0$ is a $\mathcal{T}$-relatively generating set for
a $E$-closed set $\mathcal{T}_0$ then the following conditions are
equivalent:

$(1)$ $\mathcal{T}'_0$ is the $\mathcal{T}$-least generating set
for $\mathcal{T}_0$;

$(2)$ $\mathcal{T}'_0$ is a $\mathcal{T}$-minimal generating set
for $\mathcal{T}_0$;

$(3)$ any theory in $\mathcal{T}'_0\setminus\mathcal{T}$ is
isolated by some set $(\mathcal{T}'_0\cup\mathcal{T})_\varphi$;

$(4)$ any theory in $\mathcal{T}'_0\setminus\mathcal{T}$ is
isolated by some set $(\mathcal{T}_0\cup\mathcal{T})_\varphi$;

$(5)$ any theory in $\mathcal{T}'_0\setminus\mathcal{T}$ is
isolated by some set $(\mathcal{T}'_0)_\varphi$;

$(6)$ any theory in $\mathcal{T}'_0\setminus\mathcal{T}$ is
isolated by some set $(\mathcal{T}_0)_\varphi$.}

\medskip
{\bf\em Proof.} $(1)\Rightarrow(2)$ and $(4)\Rightarrow(3)$ are
obvious.

$(2)\Rightarrow(1)$. Assume that $\mathcal{T}'_0$ is
$\mathcal{T}$-minimal but not $\mathcal{T}$-least. Then there is a
$\mathcal{T}$-relatively generating set $\mathcal{T}''_0$ such
that
$\mathcal{T}'_0\setminus(\mathcal{T}''_0\cup\mathcal{T})\ne\varnothing$
and
$\mathcal{T}''_0\setminus(\mathcal{T}'_0\cup\mathcal{T})\ne\varnothing$.
Take
$T\in\mathcal{T}'_0\setminus(\mathcal{T}''_0\cup\mathcal{T})$.

We assert that $T\in{\rm
Cl}_{E}(\mathcal{T}'_0\setminus(\{T\}\cup\mathcal{T}))$, i.e., $T$
is an accumulation point of
$\mathcal{T}'_0\setminus(\{T\}\cup\mathcal{T})$. Indeed, since
$\mathcal{T}''_0\setminus(\mathcal{T}'_0\cup\mathcal{T})\ne\varnothing$
and $\mathcal{T}''_0\subset{\rm
Cl}_E(\mathcal{T}'_0\cup\mathcal{T})={\rm
Cl}_E(\mathcal{T}'_0\setminus\mathcal{T})\cup\mathcal{T}$ (using
that $\mathcal{T}$ is $E$-closed), then by \cite[Proposition 1,
(3)]{cl} (that every finite set
$\mathcal{T}\subset\overline{\mathcal{T}}$ is $E$-closed),
$\mathcal{T}'_0\setminus\mathcal{T}$ is infinite and by
Proposition 1.1 it suffices to prove that for any $\varphi\in T$,
$((\mathcal{T}'_0\setminus(\{T\}\cup\mathcal{T}))_\varphi$ is
infinite. Assume on contrary that for some $\varphi\in T$,
$((\mathcal{T}'_0\setminus(\{T\}\cup\mathcal{T}))_\varphi$ is
finite. Then $(\mathcal{T}'_0\setminus\mathcal{T})_\varphi$ is
finite and, moreover, as $\mathcal{T}'_0$ is
$\mathcal{T}$-relatively generating for $\mathcal{T}_0$, by
Proposition 1.1, $(\mathcal{T}_0\setminus\mathcal{T})_\varphi$ is
finite, too. So $(\mathcal{T}''_0\setminus\mathcal{T})_\varphi$ is
finite and, again by Proposition 1.1, $T$ does not belong to ${\rm
Cl}_E(\mathcal{T}''_0\cup\mathcal{T})$ contradicting to ${\rm
Cl}_E(\mathcal{T}''_0)=\mathcal{T}_0$.

Since $T\in{\rm
Cl}_E(\mathcal{T}'_0\setminus(\{T\}\cup\mathcal{T}))$ and
$\mathcal{T}'_0$ is generating for $\mathcal{T}_0$, then
$\mathcal{T}'_0\setminus\{T\}$ is also generating for
$\mathcal{T}_0$ contradicting the $\mathcal{T}$-minimality of
$\mathcal{T}'_0$.

$(2)\Rightarrow(3)$. If $\mathcal{T}'_0\setminus\mathcal{T}$ is
finite then by Proposition 2.1 (3),
$\mathcal{T}'_0\setminus\mathcal{T}=\mathcal{T}_0\setminus\mathcal{T}$.
Since $\mathcal{T}_0\setminus\mathcal{T}$ is finite and
$\mathcal{T}$ is $E$-closed then for any $T\in
\mathcal{T}_0\setminus\mathcal{T}$ there is a formula $\varphi\in
T$ negating all theories in
$(\mathcal{T}_0\setminus\{T\})\cup\mathcal{T})$. Therefore,
$(\mathcal{T}_0\cup\mathcal{T})_\varphi=(\mathcal{T}'_0\cup\mathcal{T})_\varphi$
is a singleton containing $T$ and thus,
$(\mathcal{T}'_0\cup\mathcal{T})_\varphi$ isolates $T$.

Now let $\mathcal{T}'_0\setminus\mathcal{T}$ be infinite. Assume
that some $T\in\mathcal{T}'_0\setminus\mathcal{T}$ is not isolated
by the sets $(\mathcal{T}'_0\cup\mathcal{T})_\varphi$. It implies
that for any $\varphi\in T$,
$((\mathcal{T}'_0\setminus\{T\})\cup\mathcal{T})_\varphi$ is
infinite. Using Proposition 1.1 and the condition that
$\mathcal{T}$ is $E$-closed we obtain $T\in{\rm
Cl}_{E,\mathcal{T}}(\mathcal{T}'_0\setminus\{T\})$ contradicting
the $\mathcal{T}$-minimality of $\mathcal{T}'_0$.

$(3)\Rightarrow(2)$. Assume that any theory $T$ in
$\mathcal{T}'_0\setminus\mathcal{T}$ is isolated by some set
$(\mathcal{T}'_0\cup\mathcal{T})_\varphi$. By Proposition 1.1 it
implies that $T\notin{\rm
Cl}_E((\mathcal{T}'_0\setminus\{T\})\cup\mathcal{T})$. Thus,
$\mathcal{T}'_0$ is a $\mathcal{T}$-minimal generating set for
$\mathcal{T}_0$.

$(3)\Rightarrow(4)$ is obvious for finite
$\mathcal{T}'_0\setminus\mathcal{T}$. If
$\mathcal{T}'_0\setminus\mathcal{T}$ is infinite and any theory
$T$ in $\mathcal{T}'_0\setminus\mathcal{T}$ is isolated by some
set $(\mathcal{T}'_0\cup\mathcal{T})_\varphi$ then $T$ is isolated
by the set $(\mathcal{T}_0\cup\mathcal{T})_\varphi$, since
otherwise using Proposition 1.1 and the properties that
$\mathcal{T}$ is $E$-closed and $\mathcal{T}'_0$ generates
$\mathcal{T}_0$, there are infinitely many theories in
$\mathcal{T}'_0$ containing $\varphi$ that contradicts the
equality $|(\mathcal{T}'_0\cup\mathcal{T})_\varphi|=1$.

$(3)\Leftrightarrow(5)$ and $(4)\Leftrightarrow(6)$ are equivalent
since $\mathcal{T}$ is $E$-closed.~$\Box$

\medskip
{\bf Corollary 4.3.} {\em If $\mathcal{T}_j$, $j\in J$, are
$\Sigma_0$-disjoint families then $\bigcup\limits_{j\in
J}\mathcal{T}_j$ has a $\mathcal{T}_{\Sigma_0}$-least generating
set if and only if each $\mathcal{T}_j$ has a
$\mathcal{T}_{\Sigma_0}$-least generating set. Moreover, if
$\bigcup\limits_{j\in J}\mathcal{T}_j$ has a
$\mathcal{T}_{\Sigma_0}$-least generating set $\mathcal{T}_0$ then
$\mathcal{T}_0\setminus\mathcal{T}_{\Sigma_0}$ can be represented
as a disjoint union of $\mathcal{T}_{\Sigma_0}$-least generating
sets for $\mathcal{T}_j$.}

\medskip
{\bf\em Proof.} Using Theorem 4.2 it suffices to note that
$\mathcal{T}_{\Sigma_0}$ is $E$-closed and having
$\mathcal{T}_0\setminus\mathcal{T}_{\Sigma_0}$ it consists of
isolated points each of which is related to exactly one set
$\mathcal{T}_j$.~$\Box$

\medskip
Clearly, any subset of $\mathcal{T}$-least generating set is again
a $\mathcal{T}$-least generating set (for its $E$-closure). At the
same time the property ``to be a $\mathcal{T}$-least generating
set'' is preserved under finite extensions of generating sets
$\mathcal{T}'_0$ disjoint with ${\rm Cl}_E(\mathcal{T}'_0)$:

\medskip
{\bf Proposition 4.4.} {\em If $\mathcal{T}$ is a $E$-closed set,
$\mathcal{T}'_0$ is a $\mathcal{T}$-relatively generating set for
a $E$-closed set $\mathcal{T}_0$, and $\mathcal{T}_f$ is a finite
subset of $\overline{\mathcal{T}}$ disjoint with $\mathcal{T}_0$
then the following conditions are equivalent:

$(1)$ $\mathcal{T}'_0$ is the $\mathcal{T}$-least generating set
for $\mathcal{T}_0$;

$(2)$ $\mathcal{T}'_0\cup(\mathcal{T}_f\setminus\mathcal{T}_0)$ is
the $\mathcal{T}$-least generating set for the $E$-closed set
$\mathcal{T}_0\cup\mathcal{T}_f$.}

\medskip
{\bf\em Proof.} $(1)\Rightarrow(2)$. If $\mathcal{T}'_0$ is a
$\mathcal{T}$-least generating set for $\mathcal{T}_0$ then by
Theorem 4.2 each theory $T$ in $\mathcal{T}'_0\setminus
\mathcal{T}$ is isolated by some formula $\varphi_T$. Since
$\mathcal{T}_f$ is finite then each theory $T$ in
$(\mathcal{T}'_0\cup(\mathcal{T}_f\setminus\mathcal{T}_0))\setminus
\mathcal{T}$ is isolated by some formula $\psi_T$. Again by
Theorem 4.2, $\mathcal{T}'_0\cup\mathcal{T}_f$ is the
$\mathcal{T}$-least generating set for
$\mathcal{T}_0\cup\mathcal{T}_f$ which is $E$-closed in view of
Theorem 1.2.

$(2)\Rightarrow(1)$ is obvious.~$\Box$

\medskip
{\bf Theorem 4.5.} (Decomposition Theorem) {\em For any $E$-closed
sets $\mathcal{T}$ and $\mathcal{T}'$ of a language $\Sigma$ there
is a $\mathcal{T}$-relatively generating set
$\mathcal{T}'_0\cup\mathcal{T}'_1$ for $\mathcal{T}'$, which is
disjoint with $\mathcal{T}$ and satisfies the following
conditions:

$(1)$ $|\mathcal{T}'_0\cup\mathcal{T}'_1|\leq{\rm
max}\{|\Sigma|,\omega\}$;

$(2)$ $\mathcal{T}'_0$ is the least generating set for its
$E$-closure ${\rm Cl}_E(\mathcal{T}'_0)$;

$(3)$ ${\rm Cl}_E(\mathcal{T}'_0)\cap\mathcal{T}'_1=\varnothing$;

$(4)$ $\mathcal{T}'_1$ is either empty or infinite and does not
have infinite subsets satisfying $(2)$.}

\medskip
{\bf\em Proof.} We denote by $\mathcal{T}'_0$ the set of isolated
points in $\mathcal{T}'\setminus\mathcal{T}$ and by
$\mathcal{T}'_1$ the subset of
$\mathcal{T}'\setminus(\mathcal{T}\cup{\rm Cl}_E(\mathcal{T}'_0))$
with a cardinality $\leq{\rm max}\{|\Sigma|,\omega\}$ such that
each sentence belonging to a theory in
$\mathcal{T}'\setminus(\mathcal{T}\cup{\rm Cl}_E(\mathcal{T}'_0))$
belongs to a theory in $\mathcal{T}'_1$. Note that
$|\mathcal{T}'_0|$ is bounded by the number of sentences in the
language $\Sigma$, i.~e., $|\mathcal{T}'_0|\leq {\rm
max}\{|\Sigma|,\omega\}$, too. Thus the condition (1) holds and
$\mathcal{T}'_0\cup\mathcal{T}'_1$ is a $\mathcal{T}$-relatively
generating set for $\mathcal{T}'$ in view of Proposition 1.1.

By Theorem 4.2, $\mathcal{T}'_0$ is the least generating set for
${\rm Cl}_E(\mathcal{T}'_0)$. Therefore the condition (2) holds.
Now (3) and (4) are satisfied since $\mathcal{T}'_1$ is separated
from ${\rm Cl}_E(\mathcal{T}'_0)$ and does not have isolated
points.~$\Box$

\medskip
{\bf Theorem 4.6.} {\em If $T$ is a $E$-combination of some
theories $T_i$, $i\in I$, $\mathcal{T}$ is a $E$-closed set of
theories, and $|e_\mathcal{T}$-${\rm Sp}(T)|<2^\omega$, then ${\rm
Cl}_E(\mathcal{T}\cup\{T_i\mid i\in I\})$ has the
$\mathcal{T}$-least generating set.}

\medskip
{\bf\em Proof.} By Theorem 4.2 we have to show that
$\mathcal{T}'\rightleftharpoons\{T_i\mid i\in
I\}\setminus\mathcal{T}$ has a generating set, modulo
$\mathcal{T}$, of theories $T_i$ being isolated points. Assume the
contrary. Then we have sets $\mathcal{T}'_0$ and $\mathcal{T}'_1$
in terms of Theorem 4.5, where
$|\mathcal{T}'_0\cup\mathcal{T}'_1|\leq{\rm
max}\{|\Sigma|,\omega\}$ and $\mathcal{T}'_1$ is infinite. Thus
$T$ has a model $\mathcal{M}$ whose all $E$-classes satisfy
theories in $\mathcal{T}'_0\cup\mathcal{T}'_1$.

Then we can construct a $2$-tree \cite{Spr} of sentences
$\varphi_\delta$, where $\delta$ are $\{0,1\}$-tuples,
$\{\varphi_{\delta\hat{\,}0},\varphi_{\delta\hat{\,}1}\}$ are
inconsistent and
$\varphi_\delta\equiv\varphi_{\delta\hat{\,}0},\varphi_{\delta\hat{\,}1}$,
such that all $(\mathcal{T}'_1)_{\varphi_\delta}$ are infinite.
Moreover, taking negations of formulas isolating theories in
$\mathcal{T}'_1$ and applying Proposition 1.1 we can assume that
for each $f\in 2^\omega$ the sequence of formulas
$\varphi_{\langle f(0),\ldots,f(n)\rangle}$, $n\in\omega$, is
contained in a theory belonging ${\rm Cl}_E(\mathcal{T}'_1)$. Thus
$|{\rm Cl}_E(\mathcal{T}'_1)|\geq 2^\omega$ producing, by
$\mathcal{M}$, $|e_\mathcal{T}$-${\rm Sp}(T)|\geq 2^\omega$ that
contradicts the assumption $|e_\mathcal{T}$-${\rm
Sp}(T)|<2^\omega$.~$\Box$

\medskip
The following example shows that, in Theorem 4.6, the conditions
$|e_\mathcal{T}$-${\rm Sp}(T)|<2^\omega$ and the existence of the
$\mathcal{T}$-least generating set are not equivalent.

\medskip
{\bf Example 4.7.} Let $\Sigma$ be a language with predicates
$P_i$, $Q_j$, $i,j\in\omega$, of same arity (it suffices to take
the arity $0$). Now we consider a countable set of language
uniform theories $T_i$ \cite{lut} such that unique $P_i$ is
satisfied and $Q_j$ are satisfied independently for the set
$\mathcal{T}=\{T_i\mid i\in\omega\}$.

All theories $T_i$ are isolated in ${\rm Cl}_E(\mathcal{T})$ by
the formulas $\exists\bar{x}P_i(\bar{x})$. Hence, $\mathcal{T}$ is
the least generating set for ${\rm Cl}_E(\mathcal{T})$. At the
same time $|{\rm Cl}_E(\mathcal{T})|=2^\omega$ witnessed by
theories with empty predicates $P_i$ and independently satisfying
$Q_j$. Thus $|e_\mathcal{T}$-${\rm Sp}(T)|=2^\omega$ for the
theory $T$ being the $E$-combination of $T_i$,
$i\in\omega$.~$\Box$

\end{document}